\documentclass[11pt]{article}
\usepackage{amsmath, amsfonts, amssymb}
\usepackage{graphicx, amsthm}
\usepackage{dsfont}

\providecommand{\U}[1]{\protect \rule{.1in}{.1in}}

\newtheorem{theorem}{Theorem}
\newtheorem{corollary}[theorem]{Corollary}
\newtheorem{lemma}[theorem]{Lemma}
\newtheorem{proposition}[theorem]{Proposition}
\theoremstyle{definition}
\newtheorem{definition}[theorem]{Definition}

\begin{document}

\begin{center}
{\LARGE Linear structure of functions with maximal Clarke subdifferential}
\end{center}

\medskip

\begin{center}
{\large \textsc{Aris Daniilidis, Gonzalo Flores}}

\bigskip
\end{center}

\noindent \textbf{Abstract.} It is hereby established that the set of Lipschitz
functions $f:\mathcal{U}\rightarrow \mathbb{R}$ ($\mathcal{U}$ nonempty open
subset of $\ell_{d}^{1}$) with maximal Clarke subdifferential contains a
linear subspace of uncountable dimension (in particular, an isometric copy of
$\ell^{\infty}(\mathbb{N})$). This result goes in the line of a previous result of
Borwein-Wang (\cite{BW2000}, \cite{BW2003}). However, while the latter was
based on Baire category theorem, our current approach is constructive and is
not linked to the uniform convergence. In particular we establish lineability
(and spaceability for the Lipschitz norm) of the above set inside the set of
all Lipschitz continuous functions.

\vspace{0.55cm}


\vspace{0.55cm}

\textit{Secondary} 46G05, 46B04

\section{Introduction}

Let $X$ be a separable Banach space and $\mathcal{U}$ a nonempty open subset
of $X$. We denote by $\overline{B}_{\ast}$ the closed unit ball of the dual
space $X^{\ast}$ and by $||f||_{\mathrm{Lip}}$ the Lipschitz constant of a
Lipschitz function $f:\mathcal{U}\rightarrow \mathbb{R}$ (see (\ref{LipConst})
below). We also denote by $\mathrm{Lip}^{[k]}(\mathcal{U})$ the set of
Lipschitz functions $f$ defined on $\mathcal{U}$ of Lipschitz constant
$||f||_{\mathrm{Lip}}\leq k$. This space, when endowed with the metric of
uniform convergence over bounded subsets of $\mathcal{U}$, is complete.

\smallskip

In the above setting J. Borwein and X. Wang have shown in \cite{BW2000}--\cite{BW2003}, that the set of Lipschitz functions with maximal Clarke
subdifferential (that is, $\partial^{\circ}f(x)\equiv||f||_{\mathrm{Lip}}\,
\overline{B}_{\ast}$ for all $x\in \mathcal{U}$) is generic in $\mathrm{Lip}
^{[k]}(\mathcal{U})$. The result has been obtained via a standard application
of Baire's category theorem. However, this result highly depends on the chosen
metric, the reason being that wild functions with oscillating derivatives can
be obtained as uniform limits of well-behaved ones (piecewise linear or
quadratic). An explicit construction of such a wild function with maximal
Clarke subdifferential is given in \cite{BS2010}.

\smallskip

Therefore, in some generic sense, most Lipschitz functions are
\textit{Clarke-saturated} (see forthcoming Definition \ref{def-sat}), but this
genericity is strongly related to the chosen topology. To illustrate further
this fact, let us fix a nonempty compact subset $K$ of $\mathcal{U}$ and let
us consider $\mathrm{Lip}^{[k]}(K)$ as a closed subset of the Banach space
$(\mathcal{C}(K),||\cdot||_{\infty})$ (a uniform limit of Lipschitz continuous functions
of Lipschitz constant bounded by $k$ is Lipschitz). Then $||\cdot||_{\infty}$-limits of
piecewise polynomial functions in $\mathrm{Lip}^{[k]}(K)$ may give rise to
Lipschitz functions with maximal Clarke subdifferentials. A completely
different behaviour appears if one uses instead, the Lipschitz norm (see
(\ref{LipConst}) below) to describe convergence: in this case $||\cdot
||_{\mathrm{Lip}}$-limits of (piecewise) polynomials are (piecewise)
$\mathcal{C}^{1}$-functions (therefore $\partial^{\circ}f(x)\equiv \{df(x)\}$,
for all $x\in K$). The reason is that for smooth functions the Lipschitz norm
$||\cdot||_{\mathrm{Lip}}$ coincides with the norm of uniform convergence of
the derivatives and under this norm $\mathcal{C}^{1}(K)$ is a Banach subspace
of $\mathrm{Lip}(K)$.

\smallskip

If $X=\mathbb{R}^{d},$ then important subclasses of Lipschitz functions, such
as semialgebraic (more generally, o-minimal) Lipschitz functions or finite
selections of $\mathcal{C}^{d}$-smooth functions have small Clarke
subdifferentials: indeed, the aforementioned classes satisfy a Morse-Sard
theorem for their generalized critical values, see \cite[Corollary
5(ii)]{BDLS} and \cite[Theorem 5]{BDD} respectively, while every point (and
consequently every value) of a Clarke-saturated Lipschitz function is critical.

\smallskip

In this work we complement the results \cite{BS2010}, \cite{BW2000},
\cite{BW2003} by establishing a topology-independent result (Theorem~\ref{MainResult}(i)), namely, that the set of Clarke-saturated Lipschitz
functions contains an infinite dimensional linear space of uncountable
dimension; in particular it is \textit{lineable}, according to the terminology
of \cite{GQ2004}, and consequently algebraically large. Moreover,
surprisingly, $(\mathrm{Lip}(K),||\cdot||_{\mathrm{Lip}})$ contains a closed
non-separable subspace of Clarke-saturated functions, hence this set is also
\textit{spaceable. }We refer to \cite{ABPS-book} for related terminology and an
exposition on the state of the art of this trend, nowadays known as
lineability and spaceability. We also refer to \cite{ABMP2015}, \cite{M2018},
and the expository paper \cite{BPS2014}, for recent results. In some sense, our results 
have been anticipated in \cite[page 114]{BPS2014}.

\section{Preliminaries, Notation}

For any integer $d\geq1$ and real $p\in \lbrack1,\infty]$, we denote by
$\ell_{d}^{p}$ the finite-dimensional vector space $\mathbb{R}^{d}$ endowed
with the classical $p$-norm. It is a well known fact that this space is
reflexive, with $(\ell_{d}^{p})^{\ast}=\ell_{d}^{q}$, where $\frac{1}{p}
+\frac{1}{q}=1$. We denote by $\langle \cdot,\cdot \rangle:\ell_{d}^{p}
\times \ell_{d}^{q}\rightarrow \mathbb{R}$ this duality mapping .When no
confusion occures, we will simply denote the norm of $\ell_{d}^{p}$ by
$\Vert \cdot \Vert$ and the norm of $\ell_{d}^{q}$ by $\Vert \cdot \Vert_{\ast}$
(dual norm).

\smallskip

We denote by $\mathrm{Lip}(\mathcal{U})$, for $\mathcal{U}\subseteq \ell
_{d}^{p}$, the vector space of all Lipschitz functions $f:\mathcal{U}
\rightarrow \mathbb{R}$ that is, those functions for which there exists a
constant $L>0$ such that%
\begin{equation}
|f(x)-f(y)|\leq L\Vert x-y\Vert,\quad \text{for all }x,y\in \mathcal{U}.
\label{Lip}%
\end{equation}
We denote by $||f||_{\mathrm{Lip}}$ the infimum of the above constants, that
is:%
\begin{equation}
||f||_{\mathrm{Lip}}=\inf \{L>0\,:\,|f(x)-f(y)|\leq L\Vert x-y\Vert,\,
\text{for all }x,y\in \mathcal{U}\} \label{LipConst}
\end{equation}
which in turns is equivalent to
\begin{equation}
||f||_{\mathrm{Lip}}=\sup_{x,y\in \mathcal{U}\,,\,x\neq y}\frac{|f(x)-f(y)|}
{\Vert x-y\Vert}\,. \label{LipConst2}
\end{equation}
It is well-known that $||\cdot||_{\mathrm{Lip}}$ is a seminorm on
$\mathrm{Lip}(\mathcal{U})$. Fixing $x_{0}\in \mathcal{U}$ and considering the
space $\mathrm{Lip}_{x_{0}}(\mathcal{U})$ of all Lipschitz functions such that
$f(x_{0})=0$, then the aforementioned seminorm becomes a norm, and
$\mathrm{Lip}_{x_{0}}(\mathcal{U})$ a Banach space under $||\cdot
||_{\mathrm{Lip}}$.

\smallskip

Recall that every Lipschitz funcion is diferentiable almost everywhere
(Rademacher theorem). If $\mathcal{D}_{f}$ stands for the set of points where
$f$ is differentiable, and $Df(x)$ for the derivative of $f$ at $x\in
\mathcal{D}_{f}$, then the Clarke subdifferential of $f$ at $x\in \mathcal{U}$
is given by (\cite[Chapter 2]{Clarke}):
\begin{equation}
\partial^{\circ}f(x)=\overline{\mathrm{co}}\  \left \{  \lim_{x_{n}\rightarrow
x}\ Df(x_{n})\,:\, \{x_{n}\} \subseteq \mathcal{D}_{f}\right \}  .\label{a0}
\end{equation}
It follows that $\partial^{\circ}f(x)$ is a nonempty convex compact subset of
$\ell_{d}^{q}$ and for every $x^{\ast}\in \partial^{\circ}f(x)$ it holds
$\Vert x^{\ast}\Vert_{\ast}\leq||f||_{\mathrm{Lip}}$. Therefore $\partial
^{\circ}f(x)\subset||f||_{\mathrm{Lip}}\overline{B}_{\ast}$.

\begin{definition}
[Clarke-saturated function]\label{def-sat}We say that $f\in \mathrm{Lip}
(\mathcal{U})$ has a maximal Clarke subdifferential at $x_{0}\in \mathcal{U}$
whenever $\partial^{\circ}f(x_{0})\equiv||f||_{\mathrm{Lip}}\overline{B}
_{\ast}$, that is, the Clarke subdifferential equals to the closed ball of
$\ell_{d}^{q}$ centered at $0$ and with radius $||f||_{\mathrm{Lip}}$. If this
is valid for every $x\in \mathcal{U}$, we say that $f$ is Clarke-saturated.
\end{definition}

The first example of a Clarke-saturated Lipschitz function in one-dimension
has been given (up to obvious modifications) by G. Lebourg in
\cite[Proposition 1.9]{Lebourg}. The function was given by an explicit formula
based on a \textit{splitting} subset $A$ of $\mathbb{R}$ with respect to the
family of nontrivial intervals of $\mathbb{R}$, that is, a measurable subset
$A$ satisfying
\begin{equation}
0<\lambda(A\cap I)<\lambda(I),\qquad \text{for every (nontrivial) interval
}I\subset \mathbb{R}\text{,} \label{split}
\end{equation}
where $\lambda$ denotes the Lebesgue measure. An explicit construction of such
a splitting set can be found in \cite{Kirk} in a general setting (atomless
measure space). In the next section we shall enhance this construction to the
particular case of a real line and come up with a countable family of disjoint
spitting sets. This family will be paramount for the proof of our main result.

\smallskip

Let us recall that $L^{\infty}(\mathcal{U};\ell_{d}^{q})=L^{1}(\mathcal{U}
;\ell_{d}^{p})^{\ast}$ (see \cite[p. 98]{Diestel} e.g.) We shall need the
following recent result about the space $\mathrm{Lip}_{x_{0}}(\mathcal{U})$
which relates this space to some subspace of $L^{\infty}(\mathcal{U};\ell
_{d}^{q})$, the space of essentially bounded
Lebesgue-measurable functions $g:\mathcal{U}\subseteq \ell_{d}^{p}
\rightarrow \ell_{d}^{q}$. This result has been established independently in
\cite{Flores-Master} (see also \cite{F2017}) and in \cite{CK2017}.

\begin{theorem}
[isometric injection of $\mathrm{Lip}_{x_{0}}(\mathcal{U})$ into $L^{\infty
}(\mathcal{U},\ell_{d}^{q})$]\label{Isom} Let $\mathcal{U}\subseteq \ell
_{d}^{p}$ be a nonempty open convex set and $x_{0}\in \mathcal{U}$. Then, the
linear operator
\[
\left \{
\begin{array}
[c]{c}
\mathrm{\hat{D}}:\mathrm{Lip}_{x_{0}}(\mathcal{U})\rightarrow L^{\infty
}(\mathcal{U};\ell_{d}^{q})\\
\\
\mathrm{\hat{D}}f=Df\quad \text{a.e.}
\end{array}
\right.
\]
defines an isometry between \textup{$\mathrm{Lip}$}$_{x_{0}}(\mathcal{U})$ and
the following subspace of $L^{\infty}(\mathcal{U};\ell_{d}^{q})$:
\begin{equation}
\mathrm{\hat{D}}(\mathrm{Lip}_{x_{0}}(\mathcal{U}))=\left \{  g\in L^{\infty
}(\mathcal{U};\ell_{d}^{q})\,:\, \partial_{i}g_{j}=\partial_{j}g_{i}
\, \, \text{for every }i,j\in \{1,\ldots,n\} \right \}  .\label{aris}
\end{equation}

\end{theorem}

Here, $\partial_{i}g_{j}$ stands for the partial derivative of the $j$-th
component of $g$ with respect to $x_{i}$ in the sense of distributions. That
is, if $\mathcal{C}_{0}^{\infty}(\mathcal{U})$ denotes the space of test
functions (compactly supported $\mathcal{C}^{\infty}$-functions on
$\mathcal{U}$) then \eqref{aris} becomes:
\[
\int_{\mathcal{U}}g_{j}(x)\frac{\partial \varphi}{\partial x_{i}}
(x)dx=\int_{\mathcal{U}}g_{i}(x)\frac{\partial \varphi}{\partial x_{j}}
dx,\quad \text{for every }\varphi \in \mathcal{C}_{0}^{\infty}(\mathcal{U}).
\]

\section{Main result}

In this section we establish our main result which consists in exhibiting a
linear space of uncountable dimension of Clarke-saturated Lipschitz functions,
whenever $\mathcal{U}\subseteq \ell_{d}^{1}$ is a nonempty open convex set.
More precisely, endowing \textup{$\mathrm{Lip}$}$_{x_{0}}(\mathcal{U})$ with
the Lipschitz norm $||\cdot||_{\mathrm{Lip}}$ we obtain a closed subspace of
Clarke-saturated elements, which in turn implies the result thanks to Baire
theorem. Our technique is as follows: we will first prove the result for the
1-dimensional case and then we extend the construction for the $d$-dimensional
case. In both cases, we first obtain countably many linearly independent
Clarke-saturated functions in \textup{$\mathrm{Lip}$}$_{x_{0}}(\mathcal{U})$
and in the final subsection we use these functions to obtain the final result.

\subsection{The case $d=1$}

The construction for the aforementioned family of functions relies on some
basic results concerning Lebesgue measure. We refer to \cite{Fremlin} for
prerequisites in measure theory. Let us start with a typical example of a
subset of $[0,1]\ $which is closed, nowhere dense and has positive measure.

\begin{definition}
[Smith-Volterra-Cantor set]Consider the subsets $F_{n}\subset \lbrack0,1]$
defined as follows:

\begin{itemize}
\item $F_{0}=[0,1]$

\item $F_{n}$ is obtained by removing the middle open interval of length
$\frac{1}{2\cdot4^{n}}$ from each of the $2^{n}$ closed intervals whose union
is $F_{n}$.
\end{itemize}

Let $F=\bigcap_{n\geq0} F_{n}$. Then $F$ is closed and contains no intervals.
Moreover, $F$ is Lebesgue measurable with measure $1/2$.
\end{definition}

In what follows we shall use the term \textit{fat Cantor set} for any
Cantor-type set (that is, a set built in this way) with positive measure. It
is clear that this procedure can be carried out over any (open or closed)
interval, thanks to the homogeneity and invariance of the Lebesgue measure.

Let us now give the following definition:

\begin{definition}
[everywhere positive-measured set]\label{def-everywhere}A subset $A$ of
$\mathbb{R}$ is called everywhere positive-measured, if it intersects any
nontrivial interval in a set of positive measure.
\end{definition}

Notice that a set $A\ $has the splitting property (\ref{split}) for the family
of intervals of $\mathbb{R}$ if both $A$ and $\mathbb{R}\diagdown A$ are
everywhere positive-measured. The following lemma asserts the existence of a
countable partition of $\mathbb{R}$ into splitting sets.

\begin{lemma}
[countable splitting partition]\label{CounPart} There exists a countable
partition $\{A_{k}\}_{k\in \mathbb{N}}$ of $\mathbb{R}$ each of which splits
the family of intervals.
\end{lemma}

\noindent \textit{Proof.} Let us first notice that it suffices to obtain a
partition of $[0,1)$ with the above property, since we can translate those
sets over every interval of the form $[m,m+1)$, $m\in \mathbb{Z}$. To this end,
let $\{I_{n}\}_{n\in \mathbb{N}}$ be an enumeration of the subintervals of
$(0,1)$ with rational end points, say $I_{n}=(a_{n},b_{n})$. We split $I_{1}$
into two open contiguous intervals, that is, we take $c\in(a_{1},b_{1})$ and
consider the intervals $(a_{1},c)$ and $(c,b_{1})$. Then let $T_{1}^{(1)}$ and
$B^{(1)}$ be two fat Cantor sets over $(a_{1},c)$ and $(c,b_{1})$
respectively. Since $T_{1}^{(1)}\cup B^{(1)}$ is nowhere dense, there exists
$(a_{2}^{\prime},b_{2}^{\prime})\subseteq I_{2}$ such that%
\[
(a_{2}^{\prime},b_{2}^{\prime})\bigcap
\left(  T_{1}^{(1)}\cup B^{(1)}\right)  =\emptyset.
\]
We now proceed inductively as follows: Given $T_{k}^{(i)}$, $B^{(i)}$ for
$1\leq k\leq i\leq n-1$, since their union is a nowhere dense closed subset of
$(0,1)$, there exists a subinterval $(a_{n}^{\prime},b_{n}^{\prime})$ of
$I_{n}$ which is disjoint from this union. We now split the interval
$(a_{n}^{\prime},b_{n}^{\prime})$ into $n+1$ contiguous open intervals and
define $T_{k}^{(n)},B^{(n)}$ (where $k\in \{1,\ldots,n\}$) to be fat Cantor
sets over each one of these intervals. In this way we obtain inductively
disjoint fat Cantor subsets $T_{k}^{(n)}$, $B^{(n)}$ of $(0,1)$ where $1\leq
k\leq n$, and $n\in \mathbb{N}$. We then define
\[
A_{k}=\bigcup_{n\geq k}T_{k}^{(n)}\quad;\quad A_{0}=[0,1)\diagdown \left(
\bigcup_{k\geq1}T_{k}\right)  \quad;\quad B=\bigcup_{n\geq1}B^{(n)}.
\]
We claim that the family $\{A_{k}\}_{k\geq0}$ is the partition of $[0,1)$ we
are looking for.\smallskip

Indeed, the sets $\{A_{k}\}_{k\geq1}$ are mutually disjoint: Let $1\leq
k<k^{\prime}$ and assume towards a contradiction that $x\in A_{k}\cap
A_{k^{\prime}}$. Then, there exists $n\geq k$ and $n^{\prime}\geq k^{\prime}$
such that $x\in T_{k}^{(n)}$ and $x\in T_{k^{\prime}}^{(n^{\prime})}$, which
is impossible by construction. Notice further that $B\subseteq A_{0}$ (the
argument is the same as before) and that $A_{0}\subseteq \lbrack0,1)\diagdown
A_{k}$, for every $k\geq1$. Now, let $[a,b]\subseteq \lbrack0,1)$ be any
interval. For $k\geq1$, let $n\geq k$ such that $I_{n}\subseteq \lbrack a,b]$.
It follows that%
\[
\lambda(A_{k}\cap \lbrack a,b])\geq \lambda(A_{k}\cap I_{n})\geq \lambda
(T_{k}^{(n)}\cap I_{n})=\lambda(T_{k}^{(n)})>0.
\]
On the other hand
\[
\lambda(A_{0}\cap \lbrack a,b])\geq \lambda(B\cap \lbrack a,b])\geq \lambda(B\cap
I_{n})
\]%
\[
\geq \lambda(B^{(n)}\cap I_{n})=\lambda(B^{(n)})>0,
\]
yielding the result.\hfill$\square$

\bigskip

Let now $\mathcal{U}\subseteq \mathbb{R}$ be a nontrivial open interval and fix
$x_{0}\in \mathcal{U}$. Define the family of functions given by
\begin{equation}
g_{k}(x)=\mathds{1}_{A_{2k+1}}(x)-\mathds{1}_{A_{2k}}(x),\quad x\in \mathcal{U}
\label{a1}
\end{equation}
and set
\begin{equation}
f_{k}(x)=\int_{x_{0}}^{x}g_{k}(t)dt. \label{a2}
\end{equation}
We list below some properties of the family $\mathcal{F}=\{f_{k}
:k\in \mathbb{N}\}$ of functions defined by (\ref{a2}). In what follows, we
denote by $c_{00}$ the space of compactly supported sequences, that is,
$\mu=(\mu_{n})_{n\in \mathbb{N}}$ if and only if $\mathrm{supp}(\mu
):=\{n:\mu_{n}\neq0\}$ is finite.

\begin{enumerate}
\item[(i).] \label{2} $\mathcal{F}\subset \mathrm{Lip}_{x_{0}}(\mathcal{U}).$
In particular, for every $k\in \mathbb{N},$ $f_{k}$ is Lipschitz, with
$||f_{k}||_{\mathrm{Lip}}=1.$
\end{enumerate}

This is straightforward from the fact that the functions $g_{k}=f_k^{\prime}$ belong to
$L^{\infty}(\mathcal{U})$, with $\Vert g\Vert_{\infty}=1$.

\begin{enumerate}
\item[(ii).] The family $\mathcal{F}$ is linearly independent.
\end{enumerate}

Let $\mu \in c_{00}$. Then
\[
\sum_{k\in \mathbb{N}}\mu_{k}f_{k}=0\; \Longleftrightarrow \; \int_{0}^{x}\left(
\sum_{k\in \mathbb{N}}\mu_{k}g_{k}(t)\right)  dt=0,\quad \forall x\in
\mathcal{U}.
\]
In virtue of Rademacher theorem and Lebesgue differentiation theorem, the
above yields that
\[
\sum_{k\in \mathbb{N}}\mu_{k}g_{k}(x)=0,\quad \text{almost everywhere on
}\mathcal{U}.
\]
Since $\{A_{k}\}_{k\in \mathbb{N}}$ are disjoint, everywhere positive-measured
sets, we can choose $x_{k}\in A_{2k+1}\cap \mathcal{U}$, for every
$k\in \mathbb{N}$. Then $x_{k}\notin A_{2k}$, and in view of \eqref{a1} we have
$g_{k}(x_{k})=1$ and $g_{k}(x_{k^{\prime}})=0$ for $k\neq k^{\prime}$. From
this, we deduce that for every $k\in \mathbb{N}$, $\mu_{k}=0$, therefore
$\{f_{k}\}_{k\in \mathbb{N}}$ is a linearly independent family.

\begin{enumerate}
\item[(iii).] The functions $f_{k}$ are Clarke-saturated, for every
$k\in \mathbb{N}$.
\end{enumerate}

Since $f_{k}^{\prime}=g_{k}$ almost everywhere on $\mathcal{U}$, it follows
that $f_{k}^{\prime}$ takes each one of the values $\{-1,0,1\}$ on an
everywhere positive-measured (and a fortiori in a dense) subset of
$\mathcal{U}$. It follows by (\ref{a0}) that $\partial f_{k}^{\circ
}(x)=[-1,1]=\bar{B}_{\ast}(0,1)$ for every $x\in \mathcal{U}$.

\bigskip

Let us now show that the property of Clarke-saturation is inherited to linear
combinations of the family $\mathcal{F}$.

\begin{proposition}
[lineability in the 1-dim case]Every linear combination of the functions
$\{f_{k}\}_{k\in \mathbb{N}}$ has maximal Clarke sub\-diffe\-ren\-tial.
\end{proposition}

\noindent \textit{Proof.} Let $\mu \in c_{00}$ and set $f=\sum_{k\in \mathbb{N}
}\mu_{k}f_{k}$ (finite combination). Then it holds almost everywhere on
$\mathcal{U}$
\[
f^{\prime}(x)=\sum_{k\in \mathbb{N}}\mu_{k}\,f_{k}^{\prime}(x)=\sum
_{k\in \mathrm{supp}(\mu)}\mu_{k}\,g_{k}(x).
\]
Notice that for a given $x\in \mathcal{U}$ there exists at most one
$k\in \mathbb{N}$ such that $g_{k}(x)\neq0$ (namely, $g_{k}(x)=1$ or $-1$),
therefore $f^{\prime}$ can only take the values $\{ \pm \mu_{k}\}_{k\in
\mathbb{N}}$ and $0.$ Using the same argument as before, we deduce that each
of these values is taken on a dense subset of $\mathcal{U}$. Therefore
\[
\partial^{\circ}\left(  \sum_{k\in \mathbb{N}}\mu_{k}f_{k}\right)  (x)=\Vert
\mu \Vert_{\infty}[-1,1]=\bar{B}_{\ast}(0,||\mu||_{\infty}),\quad \text{for every
}x\in \mathcal{U}.
\]
Moreover,
\[
||f||_{\mathrm{Lip}}=||\sum_{k\in \mathbb{N}}\mu_{k}f_{k}||_{\mathrm{Lip}
}=\left \Vert \sum_{k\in \mathbb{N}}\mu_{k}g_{k}\right \Vert _{\infty}=\Vert
\mu \Vert_{\infty}.
\]
We conclude that this linear combination has maximal Clarke subdifferential
everywhere, that is, it is Clarke-saturated. \hfill$\square$

\subsection{The case $d>1$}

In this section we extend the above method from the 1-dimensional case to
higher dimensions. For technical reasons we equipe $\mathbb{R}^{d}$ with the
$1$-norm $||\cdot||_{1}$, so that the dual norm is $||\cdot||_{\infty}$. This
facilitates establishing Clarke-saturation. We do not know whether or not this
result remains true under a different choice of the norm. To simplify notation
we set $\ell_{d}^{1}:=(\mathbb{R}^{d},||\cdot||_{1}).$

\smallskip

Let $\mathcal{U}\subseteq \ell_{d}^{1}$ be a nonempty open convex set and let
$\mathrm{\hat{D}}$ stand for the isometry in Theorem \ref{Isom}. For
$k\in \mathbb{N}$ and $x=(x_{1},\ldots,x_{d})\in \mathcal{U}$ we define for
$k\in \mathbb{N}$ the function%
\begin{equation}
\left \{
\begin{array}
[c]{l}
G^{k}:\mathcal{U}\rightarrow \ell_{d}^{\infty}\medskip \\
G^{k}(x):=(g_{k}(x_{1}),\ldots,g_{k}(x_{d}))=(\mathds{1}_{A_{2k+1}}
(x_{1})-\mathds{1}_{A_{2k}}(x_{1}),\ldots,\mathds{1}_{A_{2k+1}}(x_{d}
)-\mathds{1}_{A_{2k}}(x_{d})).
\end{array}
\right.  \label{b1}
\end{equation}
In other words,
\[
\langle G^{k}(x),e_{i}\rangle=g_{k}(\langle x,e_{i}\rangle),
\]
where $g_{k}$ are given by \eqref{a1} and $\{e_{i}\}_{i=1,\ldots,d}$ is the
canonical basis of $\mathbb{R}^{d}$.

\smallskip

Let us first show that the functions $\{G^{k}\}_{k\in \mathbb{N}}$ are
``derivatives" of functions of \textup{$\mathrm{Lip}$}$_{x_{0}}(\mathcal{U})$.
This part relies on Theorem \ref{Isom}.

\begin{proposition}
[$G^{k}$ are derivatives]For every $k\in \mathbb{N}$, $G^{k}\in \mathrm{\hat{D}
}($\textup{$\mathrm{Lip}$}$_{x_{0}}(\mathcal{U}))$.
\end{proposition}

\noindent \textit{Proof.} Let $i,j\in \{1,\ldots,d\}$ with $i\neq j$ and
$\varphi \in \mathcal{C}_{0}^{\infty}(\mathcal{U})$. Then
\[
\int_{\mathcal{U}}\partial_{j}G_{i}^{k}(x)\, \varphi(x)\,dx=-\int_{\mathcal{U}
}G_{i}^{k}(x)\frac{\partial \varphi}{\partial x_{j}}(x)dx=-\int_{\mathcal{U}
}g_{k}(x_{i})\frac{\partial \varphi}{\partial x_{j}}(x)dx.
\]
As $\varphi \in \mathcal{C}_{0}^{\infty}(\mathcal{U})$, thanks to Fubini theorem
we can integrate first the variable $x_{j}$ and conclude that the above
integral is equal to $0$. Therefore, $\partial_{i}G_{j}^{k}=0$ whenever $i\neq
j$. In particular, $\partial_{i}G_{j}^{k}=\partial_{j}G_{i}^{k}$ in the sense
of distributions, and according to Theorem \ref{Isom} we deduce that $G^{k}
\in \mathrm{\hat{D}}($\textup{$\mathrm{Lip}$}$_{x_{0}}(\mathcal{U}))$.
\hfill$\square$

\bigskip

In view of the above proposition, we can define the family
\[
\mathcal{F}=\{f_{k}\}_{k\geq0}\subseteq \mathrm{Lip}_{x_{0}
}(\mathcal{U})
\]
as the inverse images of the family $\{G^{k}\}_{k\geq0}$, that is,
\begin{equation}
f_{k}:=\mathrm{\hat{D}}^{-1}(G^{k}),\quad \text{for every }k\in \mathbb{N}%
\text{.}\label{b2}%
\end{equation}
We now verify the same properties as in the previous section for the above functions.

\begin{itemize}
\item[(i).] $\mathcal{F}\subset \mathrm{Lip}_{x_{0}}(\mathcal{U}).$ In
particular, $||f_{k}||_{\mathrm{Lip}}=1.$
\end{itemize}

Notice that the values of $G^{k}$ are vectors $v\in \mathbb{R}^{d}$ whose
components are taking the values $\{-1,0,1\}$, each of them over everywhere
positive-measured sets. Therefore $\Vert G^{k}\Vert_{\infty}=1$ and the result
follows from the fact that $\mathrm{\hat{D}}$ is an isometry.

\begin{itemize}
\item[(ii).] The family $\mathcal{F}$ is linearly independent.
\end{itemize}

It suffices to prove that the family $\{G^{k}\}_{k\geq0}$ is linearly
independent, since $\mathrm{\hat{D}}$ is an isometry. Let $\mu \in c_{00}$
(compactly supported sequence) and assume
\[
\sum_{k\in \mathbb{N}}\mu_{k}G^{k}=0,\quad \text{that is, }\quad \sum
_{k\in \mathbb{N}}\mu_{k}G^{k}(x)=0\text{ a.e. on }\mathcal{U}.
\]
For every $k\geq0$ let $x^{k}\in(A_{k}\times \ldots \times A_{k})\cap
\mathcal{U}$. Given the definition of the functions $G^{k}$, we
have that for $i\in \{1,\ldots,d\}$
\[
\left(  \sum_{k\in \mathbb{N}}\mu_{k}G^{k}(x^{k})\right)  _{i}=\left \{
\begin{array}
[c]{cl}
\mu_{2n+1}, & \quad \text{if }k=2n+1\\
-\mu_{2n}, & \quad \text{if }k=2n.
\end{array}
\right.
\]
Since $(A_{k}\times \ldots \times A_{k})\cap \mathcal{U}$ has positive measure
everywhere, we conclude that $\mu=0$, therefore $\{G^{k}\}_{k\geq0}$ is
linearly independent and the assertion follows.

\begin{itemize}
\item[(iii).] The functions $f_{k}$ are Clarke-saturated.
\end{itemize}

Notice that every extreme point of the unit ball of $\ell_{d}^{\infty}$ is
taken as value of $G^{k}$ on a subset of $\mathcal{U}$ which has positive
measure everywhere. Since $Df_{k}=G^{k}$ almost everywhere on $\mathcal{U}$,
we conclude that $\partial^{\circ}f_{k}(x)=\bar{B}_{\ast}$, for all
$x\in \mathcal{U}$.

\bigskip

Similarly to the 1-dimensional case we now establish that Clarke-saturation is
preserved under linear combinations of elements of $\mathcal{F}$.

\begin{proposition}
[lineability]Every linear combination of the functions $(f_{k})_{k\in
\mathbb{N}}$ has maximal Clarke sub\-diffe\-ren\-tial.
\end{proposition}

\noindent \textit{Proof.} Let $\mu \in c_{00}$. Then we have
\[
D\left(  \sum_{k\in \mathbb{N}}\mu_{k}f_{k}\right)  (x)=\sum_{k\in \mathbb{N}
}\mu_{k}G^{k}(x),\quad \text{for a.e. }x\in \mathcal{U}\text{.}
\]
The values of this last function are exclusively vectors $v\in \mathbb{R}$ with
components in the set $\{ \pm \mu_{k}\,:\,k\geq0\}$. Moreover, each component
takes each one of the values $\{ \pm \mu_{k}\}_{k\in \mathbb{N}}$ on subsets of
$\mathcal{U}$ which have everywhere positive measure. It follows readily from
(\ref{a0}) that for every $x\in \mathcal{U}$
\[
\partial^{\circ}\left(  \sum_{k\in \mathbb{N}}\mu_{k}f_{k}\right)  (x)=\Vert
\mu \Vert_{\infty}B_{\ast}.
\]
In addition, using the isometry $\mathrm{\hat{D}}$ we deduce:
\[
||f||_{\mathrm{Lip}}=||\sum_{k\in \mathbb{N}}\mu_{k}f_{k}||_{\mathrm{Lip}
}=\left \Vert \sum_{k\in \mathbb{N}}\mu_{k}G^{k}\right \Vert _{\infty}=\Vert
\mu \Vert_{\infty}.
\]
The proof is complete.\hfill$\square$

\subsection{The space of Clarke-saturated functions}

In the previous section we constructed a countable family of linearly
independent Clarke-saturated functions $f_{k}$ belonging to
\textup{$\mathrm{Lip}$}$_{x_{0}}(\mathcal{U})$, where $\mathcal{U}
\subseteq \ell_{d}^{1}$ is a nonempty open convex set and $x_{0}\in \mathcal{U}
$. We shall now describe in terms of the isometry $\mathrm{\hat{D}}$ (Theorem
\ref{Isom}) the closure of the space generated by these functions. In what
follows we denote by $\ell^{\infty}(\mathbb{N})$ the (nonseparable) Banach
space of bounded sequences.

\begin{proposition}
\label{isometry} Let $T:\ell^{\infty}(\mathbb{N})\rightarrow L^{\infty
}(\mathcal{U};\ell_{d}^{\infty})$ given by
\[
T(\mu)=\sum_{k\geq0}\mu_{k}G^{k},\quad \text{for all }\mu=(\mu_{n}
)_{n\in \mathbb{N}}\in \ell^{\infty}(\mathbb{N}).
\]
Then $T$ is well defined and establishes a linear isometric injection of
$\ell^{\infty}(\mathbb{N})$ into $L^{\infty}(\mathcal{U};\ell_{d}^{\infty})$.
\end{proposition}

\noindent \textit{Proof.} Let $\{A_{k}\}_{k\in \mathbb{N}}$ be the countable
partition of $\mathbb{R}$ given by Lemma \ref{CounPart}. Let $x=(x_{1}
,\ldots,x_{d})\in \mathcal{U}$. Since each $A_{k}$ is everywhere
positive-measured, there exists $j_{1},\ldots,j_{d}\geq0$ such that $x_{i}\in
A_{j_{i}}$, for $i\in \{1,\ldots,d\}$. This implies that the sum
\[
\sum_{k\geq0}\mu_{k}G^{k}(x)
\]
is finite for every $x\in \mathcal{U}$, with norm less than or equal to
$\Vert \mu \Vert_{\infty}$. Therefore $T(\mu)\in L^{\infty}(\mathcal{U};\ell
_{d}^{\infty})$, with $\Vert T\mu \Vert_{\infty}\leq \Vert \mu \Vert_{\infty}$.
Moreover, if $x\in(A_{2n+1}\times \ldots \times A_{2n+1})$ and $x^{\prime}
\in(A_{2n}\ldots \times A_{2n})$ then
\[
T(\mu)(x)=-T(\mu)(x^{\prime})=(\mu_{k},\ldots,\mu_{k}),
\]
which leads to $\Vert T\mu \Vert_{\infty}=\Vert \mu \Vert_{\infty}$. Since $T$ is
obviously linear, it follows that $T$ is a linear isometry between
$\ell^{\infty}(\mathbb{N})$ and $T(\ell^{\infty}(\mathbb{N}))$.\hfill$\square$

\bigskip

Now, we state the relation between $T(\ell^{\infty}(\mathbb{N}))$ and
$\mathrm{\hat{D}}($\textup{$\mathrm{Lip}$}$_{x_{0}}(\mathcal{U}))$. This
relation is obtained in a similar way as in the case of linear combinations.

\begin{proposition}
\label{subset} $T(\ell^{\infty}(\mathbb{N}))\subseteq \mathrm{\hat{D}}
($\textup{$\mathrm{Lip}$}$_{x_{0}}(\mathcal{U}))$.
\end{proposition}

\noindent \textit{Proof.} Let $\mu \in \ell^{\infty}(\mathbb{N})$. We need to
prove that $T(\mu)$ is the gradient of some Lipschitz function. Let
$i,j\in \{1,\ldots,d\}$ with $i\neq j$. Then
\[
(T\mu)_{i}(x)=\sum_{k\geq0}\mu_{k}g_{k}(x_{i}).
\]
If $\varphi \in \mathcal{C}_{0}^{\infty}(\mathcal{U})$, we have that
\[
\left \langle \partial_{j}(T\mu)_{i},\varphi \right \rangle =-\int_{\mathcal{U}
}\left(  \sum_{k\geq0}\mu_{k}g_{k}(x_{i})\right)  \frac{\partial \varphi
}{\partial x_{j}}(x)dx=-\int_{\mathcal{U}}\sum_{k\geq0}\left(  \mu_{k}
g_{k}(x_{i})\frac{\partial \varphi}{\partial x_{j}}(x)\right)  dx.
\]
We define for $n\geq0$
\[
\psi_{n}(x)=\sum_{k=0}^{n}\left(  \mu_{k}g_{k}(x_{i})\frac{\partial \varphi
}{\partial x_{j}}(x)\right)  \quad \text{and}\quad \psi(x)=\sum_{k\geq0}\left(
\mu_{k}g_{k}(x_{i})\frac{\partial \varphi}{\partial x_{j}}(x)\right)  .
\]
Notice that for $x=(x_{1},\ldots,x_{d})\in \mathcal{U}$ and $i\in
\{1,\ldots,d\}$ we have
\[
g_{k}(x_{i})\neq0\; \Longleftrightarrow \;x_{i}\in A_{2k+1}\cup A_{2k}
\]
and in this case $g_{k^{\prime}}(x_{i})=0,$ for all $k^{\prime}\neq k.$
Therefore, there exists some $N\geq0$ large enough such that
\[
\psi_{n}(x)=\sum_{k=0}^{n}\mu_{k}g_{k}(x_{i})\frac{\partial \varphi}{\partial
x_{j}}(x)=\left \{
\begin{array}
[c]{cl}
0,\medskip & n<N\\
\mu_{N}\ g_{N}(x_{i})\  \frac{\partial \varphi}{\partial x_{j}}(x), & n\geq N
\end{array}
\right.
\]
yielding
\[
\psi_{n}\rightarrow \psi \text{ \ (pointwise) \quad and \quad}|\psi_{n}
|\leq||\mu||_{\infty}\left \vert \frac{\partial \varphi}{\partial x_{j}
}\right \vert \in L^{1}(\mathcal{U}).
\]
In virtue of the Lebesgue dominated convergence theorem, we have that
\[
\left \langle \partial_{j}(T\mu)_{i},\varphi \right \rangle =-\sum_{k\geq
0}\left(  \int_{\mathcal{U}}\mu_{k}g_{k}(x_{i})\frac{\partial \varphi}{\partial
x_{j}}(x)dx\right)  .
\]
But thanks to Fubini theorem, we can integrate first with respect to the
$x_{j}$ variable and since $\varphi$ has compact support, we conclude that all
the integrals are equal to $0$. Then $\partial_{j}(T\mu)_{i}=0$ whenever
$i\neq j$, which leads to $T(\mu)\in \mathrm{\hat{D}}($\textup{$\mathrm{Lip}$
}$_{x_{0}}(\mathcal{U}))$. \hfill$\square$

\bigskip

\begin{proposition}
Let $f\in$\textup{$\mathrm{Lip}$}$_{x_{0}}(\mathcal{U})$ be such that
$\mathrm{\hat{D}}f=T(\mu)$. Then $f$ is Clarke-saturated.
\end{proposition}

\noindent \textit{Proof.} It suffices to notice that
\[
\Vert f\Vert_{\mathrm{Lip}}=\Vert \mathrm{\hat{D}}f\Vert_{\infty}=\Vert
T(\mu)\Vert_{\infty}=\Vert \mu \Vert_{\infty}
\]
and that for every extreme point $v$ of the dual ball $\overline{B}_{\ast}$
and $k\geq0$ there exists an everywhere positive-measured set $A\subseteq
\mathcal{U}$ such that
\[
\mathrm{\hat{D}}f(x)=T\mu(x)=\mu_{k}v\quad \text{for every }x\in A.
\]
Since $f$ is differentiable almost everywhere, we conclude that
\[
\partial^{\circ}f(x)=\Vert \mu \Vert_{\infty}\overline{B}_{\ast}=\overline
{B}_{\ast}(0,\Vert f\Vert_{\mathrm{Lip}}).
\]
which finishes the proof.\hfill$\square$

\bigskip

We are ready to state our main result.

\begin{theorem}
[Spaceability of Clarke-saturated functions]\label{MainResult} Let $d\geq1$
and $\mathcal{U}\subseteq \ell_{d}^{1}$ be a nonempty open convex set. Then,

\begin{enumerate}
\item[(i)] {\rm (lineability)} The space \textup{$\mathrm{Lip}$}$(\mathcal{U})$ of
Lipschitz functions contains a linear subspace of Clarke-saturated functions
of uncountable dimension. 

\item[(ii)] {\rm (spaceability)} For any $x_{0}\in \mathcal{U}$, the Banach space
$($\textup{$\mathrm{Lip}$}$_{x_{0}}(\mathcal{U}),||\cdot||_{\mathrm{Lip}})$
contains a (proper) linear subspace of Clarke-saturated functions isometric to
$\ell^{\infty}(\mathbb{N})$.
\end{enumerate}
\end{theorem}

In particular, if $\mathcal{F}=\{f_{k}:k\in \mathbb{N}\}$ is the family defined
in (\ref{b2}), then $\text{span}\{f_{k}\}$ is isometric to $c_{00}$ while
$\overline{\text{span}}\{f_{k}\}$ is isometric to $c_{0}(\mathbb{N})$ (the
Banach space of null sequences).

\medskip

\noindent \textit{Proof.} Thanks to Proposition \ref{isometry} and Proposition
\ref{subset}, we deduce that $\ell^{\infty}(\mathbb{N})$ is isometric to the
subspace
\[
Z=\mathrm{\hat{D}}^{-1}(T(\ell^{\infty}(\mathbb{N})))
\]
of $\mathrm{Lip}_{x_{0}}(\mathcal{U})$. This subspace is proper (any strictly
differentiable function $h\in \textup{Lip}_{x_{0}}(\mathcal{U})\setminus \{0\}$
does not belong to $Z$). This proves (ii), and yields directly that
Clarke-saturated functions contain a linear subspace of uncountable dimension.
Therefore (i) holds, since $\mathrm{Lip}_{x_{0}}(\mathcal{U})$ is a linear subspace
of $\mathrm{Lip}(\mathcal{U})$. Finally, an easy computation shows that if
$\mu \in c_{00},$ then $\mathrm{\hat{D}}^{-1}\left(  T(\mu)\right)  \in
$span$\{f_{k}\},$ whence $c_{00}$ is isometric to $\text{span}\{f_{k}\}.$ It
follows readily by continuity that $c_{0}(\mathbb{N})$ is isometric to
$\overline{\text{span}}\{f_{k}\}$. \hfill$\square$

\bigskip

We conclude this work with the following straightforward consequence of
Theorem \ref{MainResult}.

\begin{corollary}
Let $p\in \mathbb{R}^{d}$ and $r>0$. Then, there exists $f\in$
\textup{$\mathrm{Lip}$}$(\mathcal{U})$ such that $\partial^{\circ
}f(x)=\overline{B}_{\ast}(p,r)$ for every $x\in \mathcal{U}$.
\end{corollary}

\noindent \textit{Proof.}  Let $\mu \in \ell^{\infty}$ be such that $\| \mu \|_{\infty}=r$. Set $h_{1}=D^{-1}T\mu$ and $h_{2}=\left \langle p,\cdot \right \rangle$. Then $\partial^{\circ}h_{1}(x) = \overline{B}_{*}(0,r)$ and $\partial^{\circ}h_{2}(x) = \{p\}$ for every $x\in \mathcal{U}$, where we used that $h_{2}$ is strictly differentiable. Again thanks to that, if $f=h_{1}+h_{2}$ then for every $x\in \mathcal{U}$
	\[ \partial^{\circ}f(x) = \partial^{\circ}(h_{1}+h_{2})(x) = \partial^{\circ}h_{1}(x) + \partial^{\circ}h_{2}(x) = \overline{B}_{*}(p,r). \]
The proof is complete. \hfill$\square$
\bigskip

\textbf{Acknowledgments.} This research has been supported by the following grants: CMM-AFB170001, FONDECYT 1171854 (Chile) and MTM2014-59179-C2-1-P (Spain). The second author has been supported by a doctorate scholarship of CONICYT (Chile).

\bigskip

\vspace{0.8cm}

\noindent Aris DANIILIDIS, Gonzalo FLORES

\medskip

\noindent DIM--CMM, UMI CNRS 2807\newline Beauchef 851, FCFM, Universidad de
Chile \smallskip

\noindent E-mail: \texttt{arisd@dim.uchile.cl}\, \, \,;\, \,
\texttt{gflores@dim.uchile.cl} \newline \noindent
\texttt{http://www.dim.uchile.cl/\symbol{126}arisd}

\medskip

\noindent Research supported by the grants: \newline CMM AFB170001 (Chile),
FONDECYT 1171854 (Chile), MTM2014-59179-C2-1-P (Spain).\newline

\end{document}